\documentclass[11 pt]{amsart}
\usepackage{amsmath}
\usepackage{latexsym,mathrsfs}
\usepackage{mathdots}
\usepackage{color}
\usepackage[inline]{enumitem} 
\usepackage{tikz}
\usepackage{tikz-cd}
\usepackage[T1]{fontenc}
 \usepackage[utf8]{inputenc}
 \usepackage{mathtools} 
  \usepackage{fdsymbol} 
 \usepackage{hyperref}
 \usepackage{subfiles}
\usepackage[capitalise]{cleveref}

\usepackage[margin=1.4in]{geometry}

\usetikzlibrary{positioning} 
\usetikzlibrary{decorations.pathreplacing}
\usetikzlibrary{arrows, decorations.markings} 
\pgfarrowsdeclarecombine{twotriang}{twotriang}{stealth}{stealth}%
{stealth}{stealth}
\tikzset{arrow/.style={-stealth}}
\tikzset{arrowshorter/.style={-stealth, shorten <=2pt, shorten >=2pt}}
\tikzset{arrowmuchshorter/.style={-stealth, shorten <=7pt, shorten >=6pt}}
\tikzset{mono/.style={>-stealth}} 
\tikzset{epi/.style={-twotriang}} 
\tikzset{twoarrowlonger/.style={double,double distance=1.5pt,
shorten <=5pt,shorten >=6pt,
decoration={markings,mark=at position -4pt with {\arrow[scale=1.75]{>}}},
preaction={decorate}}} 

\tikzset{twoarrow/.style={double,double distance=1.5pt,
shorten <=6pt,shorten >=7pt, 
decoration={markings,mark=at position -4pt
with {\arrow[scale=1.75]{>}}},
preaction={decorate} 
}
}
\tikzset{%
    symbol/.style={%
        draw=none,
        every to/.append style={%
            edge node={node [sloped, allow upside down, auto=false]{$#1$}}}
    }
}

\tikzset{mapstikz/.style={-stealth, 
decoration={markings,mark=at position 0pt with {\arrow[scale=0.5]{|}}}, preaction={decorate}}}
\usetikzlibrary{backgrounds,shapes,arrows,calc,patterns} 

\pgfdeclarelayer{background}
\pgfsetlayers{background,main}

\theoremstyle{definition}
\newtheorem{thm}{Theorem}[section] 
\makeatletter\let\c@thm\c@thm\makeatother
\newtheorem{cor}{Corollary}[section]
\makeatletter\let\c@cor\c@thm\makeatother
\newtheorem{lem}{Lemma}[section]
\makeatletter\let\c@lem\c@thm\makeatother
\newtheorem{prop}{Proposition}[section]
\makeatletter\let\c@prop\c@thm\makeatother

\makeatletter\let\c@claim\c@thm\makeatother

\makeatletter\let\c@conjecture\c@thm\makeatother

\newtheorem*{unnumberedtheorem}{Theorem}

\newtheorem*{unnumberedproposition}{Proposition}

\theoremstyle{definition}

\newtheorem{defn}{Definition}[section]
\makeatletter\let\c@defn\c@thm\makeatother
\newtheorem{const}{Construction}[section]
\makeatletter\let\c@const\c@thm\makeatother
\newtheorem{notn}{Notation}[section]
\makeatletter\let\c@notn\c@thm\makeatother

\makeatletter\let\c@convention\c@thm\makeatother

\newtheorem{rmk}{Remark}[section]
\makeatletter\let\c@rmk\c@thm\makeatother

\makeatletter\let\c@ex\c@thm\makeatother

\makeatletter\let\c@observation\c@thm\makeatother

\makeatletter\let\c@warning\c@thm\makeatother

\makeatletter\let\c@digression\c@thm\makeatother

\makeatletter\let\c@answ\c@thm\makeatother

\makeatletter\let\c@answ\c@thm\makeatother

\makeatletter\let\c@aside\c@thm\makeatother

\makeatletter
\let\c@equation\c@thm
\numberwithin{equation}{section}
\makeatother

\newcommand{\newrefformat}[2]{}

\crefname{lem}{Lemma}{Lemmas}
\crefname{thm}{Theorem}{Theorems}
\crefname{defn}{Definition}{Definitions}
\crefname{notn}{Notation}{Notations}
\crefname{const}{Construction}{Constructions}
\crefname{prop}{Proposition}{Propositions}
\crefname{rmk}{Remark}{Remarks}
\crefname{cor}{Corollary}{Corollaries}
\crefname{equation}{Display}{Displays}
\crefname{ex}{Example}{Examples}
\crefname{thmalph}{Theorem}{Theorems}
\crefname{answ}{Answer}{Answers}
\crefname{crzcond}{Crazy Condition}{Crazy Conditions}

\newcommand{\cC}{\mathcal{C}}

\newcommand{\cE}{\mathcal{E}}

\newcommand{\cM}{\mathcal{M}}

\newcommand{\cS}{\mathcal{S}}

\newcommand{\cat}{\cC\!\mathit{at}}
\newcommand{\set}{\cS\!\mathit{et}}
\newcommand{\sset}{\mathit{s}\set}

 \newcommand{\Fun}{\mathbf{Fun}}

\DeclareMathOperator{\tr}{tr}
\DeclareMathOperator{\core}{core}

\newcommand{\noloc}{%
  \nobreak\mskip6mu\mathord{:}\mskip2mu\relax
}

\newcommand{\adjunction}[4]{%
  \if\relax\detokenize{#1}\relax\else
    #1\nobreak\mskip2mu\mathord{:}\mskip6mu
  \fi
  #2\rightleftarrows#3%
  \if\relax\detokenize{#4}\relax\else
    \noloc#4%
  \fi
}

\DeclareMathOperator{\colim}{colim}
\DeclareMathOperator{\id}{id}

\DeclareMathOperator{\Mor}{Mor}

\author{Viktoriya Ozornova}
\address{Max Planck Institute for Mathematics, Bonn, Germany}
\email{viktoriya.ozornova@mpim-bonn.mpg.de}

\author{Martina Rovelli}
\address{
University of Ottawa, 
Ottawa,
Canada
}
\email{mrovelli@uottawa.ca}

\author{Tashi Walde}
\address{Universität Regensburg, Regensburg, Germany}
\email{tashi.walde@ur.de}

\tikzcdset{arrow style=tikz,
           diagrams={>=stealth}
           }

\title{Models of \texorpdfstring{$n$}{n}-dimensional categories as \texorpdfstring{$n\to\infty$}{n tends to infinity}}

\begin{document}

\begin{abstract}
  We verify that saturated complicial sets model
  the right limit of the $\infty$-categories of $(\infty,n)$-categories,
  whereas the canonical comparison with the left limit is not an equivalence.
  By contrast the category of strict $\omega$-categories
  is equivalent to both the right and left limits
  of the categories of strict $n$-categories.
\end{abstract}

\maketitle

\addtocontents{toc}{\protect\setcounter{tocdepth}{1}}
\tableofcontents

\section*{Introduction}

Suppose we have a sequence of inclusions of $\infty$-categories
\[
  \mathscr C_0 \hookrightarrow \mathscr C_1 \hookrightarrow \cdots \hookrightarrow \mathscr C_n
  \hookrightarrow \mathscr C_{n+1} \hookrightarrow \cdots .
\]
A natural question is how to combine these into a single $\infty$-category that reflects the behavior of the sequence as $n\to\infty$. While one can always form the colimit of the inclusions, in the examples of interest this colimit typically does not capture an interesting limiting object, and we are therefore led to consider alternative notions of limit for such sequences.

In various cases of interest the inclusions admit left or right adjoints.
This makes additional constructions available:
one can form the limit of the $\mathscr C_n$
along the left adjoints or along the right adjoints,
obtaining the \emph{left limit} $\lim_{n\geq0}^{L}\mathscr C_n$
or the \emph{right limit} $\lim_{n\geq0}^{R}\mathscr C_n$, respectively.
When all the $\mathscr C_n$ are presentable,
these constructions correspond to computing
$\colim_{n\geq 0}\mathscr C_n$ not in $\infty$-categories,
but in $\Pr^R$ or $\Pr^L$, respectively.
We discuss the framework needed to recognize when an $\infty$-category is equivalent to the left or right limit in a general setting, mimicking the formalism developed in \cite[\textsection5.5.6]{htt}.

As a first and relatively straightforward application of this general framework, we consider the sequence $\mathscr C_n=n\cat$ of categories of strict $n$-categories. We verify as \cref{LimitIndStrict,LimitCoIndStrict} that the category $\omega\cat$ of strict $\omega$-categories from \cite{StreetOrientedSimplexes} is equivalent to both the right and the left limit of the categories $n\cat$. In particular, the two limits are equivalent.

\begin{unnumberedproposition}
Let $\omega\cat$ denote the category of strict $\omega$-categories. There are equivalences
\[
\lim_{n\geq0}^{R}n\cat\quad\simeq\quad\omega\cat\quad\simeq\quad \lim_{n\geq0}^{L}n\cat.
\]
\end{unnumberedproposition}

The right limit description of strict $\omega$-categories is consistent with the framework of \cite[\textsection1.4]{LeinsterBook}, \cite[\textsection14.2.6]{PolygraphBook}, \cite[\textsection2]{JFR} and \cite[\textsection2.1]{CampionPasting},
where this description is implicitly assumed or taken as a definition.

We next turn to our main case of interest, the sequence $\mathscr C_n=\mathscr Comp_n$, where $\mathscr Comp_n$ is the $\infty$-category presented by $n$-complicial sets \cite{VeritySlides,EmilyNotes,OR,RiehlVerityBook}. These model $(\infty,n)$-categories \cite[Thm.~3.4.3.2]{LoubatonThesis}.
The right limit in this setting has been widely considered in the literature, for example in \cite[\textsection2.6]{Stefanich}, \cite[\textsection3.1]{GepnerHeine}, \cite{Goldthorpe}, and \cite[\textsection2.1]{CampionPasting}.
We verify as \cref{LimitInd} that Verity's saturated complicial sets from \cite{VeritySlides,EmilyNotes,OR,RiehlVerityBook} model the right limit of the $\infty$-categories $\mathscr Comp_n$.
As such, complicial sets provide a model for $(\infty,n)$-categories in the limit as $n\to\infty$.

\begin{unnumberedtheorem}
  \label{thm:main}
Let $\mathscr Comp$ denote the $\infty$-category of saturated complicial sets. There is an equivalence
\[\lim_{n\geq0}^{R}\mathscr Comp_n\quad\simeq\quad\mathscr Comp.\]
\end{unnumberedtheorem}

Apart from the general $\infty$-categorical techniques used to establish the formal framework,
our proof uses only a routine strictification lemma about $\mathbb{N}$-sequences in model categories,
and a small number of basic facts about the complicial models.

The identification of saturated complicial sets with the right limit is also compatible with Verity's original vision \cite{VeritySlides}.
The analogous statement for $\Theta$-spaces was proved by Loubaton \cite{LoubatonThesis},
and our result can also be deduced from his through the known equivalence of models;
see \cref{ComparisonWithFelixFramework} for more details.
An analogous statement in the semi-strict setting also appears in \cite[\textsection4.3]{henry-loubaton} (cf.~also \cite[pp.~8-9]{Loubaton4}).

We also remark that the formal framework still yields an adjunction
\[
\adjunction{}{\mathscr Comp}{\lim_{n\geq0}^{L}\mathscr Comp_n}{},
\]
which is identified with the proper reflective localization of \cite{ORW1}
(also constructed in \cite{GH-hyper}); see \cref{LimitCoind}.
The complicial set $N^{\mathrm{RS}}\cE$ considered in \cite{ORW2}
is an explicit presentation for an object in the kernel of this localization;
through the main theorem it thus legitimately presents
a nontrivial $(\infty,\infty)$-category whose finite truncations are all trivial.

\subsubsection*{Acknowledgments}
MR gratefully acknowledges support from NSF grant DMS-2203915 and an NSERC Discovery Grant. This work was partially undertaken at the Isaac Newton Institute for Mathematical Sciences, Cambridge, during the programme \emph{Equivariant homotopy theory in context} (EPSRC grant EP/Z000580/1) and at the University of Hamburg, during events under CRC~1624 \emph{Higher structures, moduli spaces, and integrability}. This work was partially supported by the CRC~1785 \emph{Generalised Motivic Methods in Geometry} (Universit\"at Regensburg, funded by the DFG). We are grateful to Thorger Geiß, Lennart Meier and Markus Spitzweck for insightful conversations on this topic.

\subsubsection*{AI use disclosure}
The mathematical content of this paper was developed and written by the authors. GPT-6 was used for mathematical discussion and proofreading. The AI assistant suggested various small local improvements and performed routine editing tasks such as formatting and fixing typos; all its edits were checked and approved by the authors.

\section{Formal framework}

\label{Framework}

\newcommand{\SecCart}[2]{\mathrm{Sec}_{#1}\mathrm{Cart}_{#2}}

\begin{const}
\label{InductiveLimit}
Given a family of adjunctions
\[
\{\adjunction{J_n}{\mathscr{C}_n}{\mathscr{C}}{R_n}\}_{n\geq 0},
\]
the \emph{right limit} $\lim_{n\geq 0}^{R}\mathscr{C}_n$ is the limit of the tower
\begin{equation}
\label{Tower}
\mathscr{C}_0 \xleftarrow{R_0J_1} \mathscr{C}_1 \xleftarrow{} \ldots \xleftarrow{} \mathscr{C}_{n-1} \xleftarrow{R_{n-1}J_n} \mathscr{C}_n \xleftarrow{R_nJ_{n+1}} \mathscr{C}_{n+1} \xleftarrow{} \ldots.
\end{equation}
\end{const}

The goal of this section is to develop some criteria that let us check
that the right limit is equivalent to the $\infty$-category $\mathscr{C}$ itself.
The constructions and proofs in this section are similar to Lurie’s discussion
of Postnikov towers in \cite[§5.5.6]{htt},
especially \cite[Remark 5.5.6.24]{htt}.
We use many of the same techniques.

\begin{defn}
A family of adjunctions $\{\adjunction{J_n}{\mathscr{C}_n}{\mathscr{C}}{R_n}\}_{n\geq 0}$
is a \emph{preinductive tower} for a cocomplete $\infty$-category $\mathscr{C}$ if:
\begin{enumerate}[label=(A\arabic*), ref=(A\arabic*)]

\item \label{A1}
For all $n\geq 0$, the functor $J_n$ is fully faithful; equivalently, for all $X$ in $\mathscr{C}_n$, the unit is an equivalence
\[
\eta_nX\colon X \xrightarrow{\simeq} R_nJ_nX.
\]

\item \label{A2}
For all $n>0$ and all $X$ in $\mathscr{C}$, the counit induces an equivalence
\[
R_{n-1}\varepsilon_nX\colon
R_{n-1}J_nR_nX \xrightarrow{\simeq} R_{n-1}X.
\]
\end{enumerate}
\end{defn}

It might be more illuminating to present the same data in a slightly different way:

\begin{prop}
\label{prop:preinductive-reflective-subcats}
By replacing each $\mathscr{C}_n$ with its essential image in $\mathscr{C}$, the data of a preinductive tower is equivalent to that of a cocomplete $\infty$-category $\mathscr{C}$ equipped with a nested sequence
\[
\mathscr{C}_0\subseteq \mathscr{C}_1\subseteq \cdots\subseteq \mathscr{C}_n\subseteq \cdots \subseteq\mathscr{C}
\]
of full subcategories, all of which are coreflective (i.e., all inclusions have a right adjoint). Under this correspondence, the right limit $\lim_{n\geq 0}^{R}\mathscr{C}_n$ is the limit of the tower
\[\mathscr{C}_0\xleftarrow{R_0^1}\mathscr{C}_1\xleftarrow{R_1^2}\cdots\xleftarrow{}\mathscr{C}_{n-1}\xleftarrow{R_{n-1}^{n}}\mathscr{C}_{n}\xleftarrow{R_n^{n+1}} \cdots\]
where each $R_n^{n+1}$ is a right adjoint to the inclusion $\mathscr{C}_n\subseteq \mathscr{C}_{n+1}$.
\end{prop}

\Cref{prop:preinductive-reflective-subcats} is a direct consequence
of a few basic lemmas about adjunctions, which we shall prove right below.

\begin{lem}
\label{lem:factor-adjunction-ff}
Consider adjunctions of $\infty$-categories
\[
\adjunction{F_1}{\mathscr A_1}{\mathscr C}{R_1}
\quad\text{and}\quad
\adjunction{F_2}{\mathscr A_2}{\mathscr C}{R_2},
\]
and assume that $F_2$ is fully faithful. The following four conditions are equivalent:
\begin{enumerate}[leftmargin=*, label=(\arabic*), ref=(\arabic*)]
\item \label{it:F-1-maps-to-image-of-F-2}
The functor $F_1$ maps each object of $\mathscr A_1$
into the essential image
$F_2(\mathscr A_2)\subseteq\mathscr C$ of $F_2$.
\item \label{it:F-1-factors}
The functor $F_1$ factors through $F_2$,
i.e., there is a functor $F'\colon \mathscr A_1\to \mathscr A_2$
with a natural equivalence $F_2F'\simeq F_1$.
\item \label{it:R-1-factors} The functor $R_1$ factors through $R_2$,
i.e., there is a functor $R'\colon \mathscr A_2\to \mathscr A_1$
with a natural equivalence $R'R_2\simeq R_1$.
\item \label{it:R-1-inverts-epsilon-2}
For each object $X$ of $\mathscr C$, the map $R_1\varepsilon_2 X\colon R_1F_2R_2 X \to R_1 X$
is an equivalence.
\end{enumerate}
If these equivalent conditions are satisfied, then there are natural equivalences
\[
F'\simeq R_2F_1 \colon \mathscr A_1\to \mathscr A_2
\quad\text{and}\quad
R'\simeq R_1F_2 \colon \mathscr A_2\to \mathscr A_1
\]
and an adjunction $F'\dashv R'$.
\end{lem}

The proof is a routine verification using the basic features of adjunctions
(natural transposition equivalence on mapping spaces, unit, counit, triangle identities);
we record it for completeness.

\begin{proof}
  \begin{itemize}[leftmargin=1cm]
  \item[]
  \item {``\ref{it:F-1-maps-to-image-of-F-2}$\implies$\ref{it:F-1-factors}'':}
    We can define $F'\coloneq R_2F_1$ together with the transformation
    $\varepsilon_2F_1\colon F_2R_2F_1 \to F_1$,
    which is an equivalence because $F_1a_1$ lies in the essential image of $F_2$ for each object $a_1$ of $\mathscr A_1$.
  \item {``\ref{it:F-1-factors}$\implies$\ref{it:R-1-factors}'':}
    Assume that we have $F'\colon \mathscr A_1\to\mathscr A_2$ with $F_2F'\simeq F_1$.
    We get an adjunction $F'\dashv R_1F_2$ from the natural equivalence
    \[
      \mathscr A_1(-,R_1F_2-)\simeq\mathscr C(F_1-,F_2-)\simeq\mathscr C(F_2F'-,F_2-)\simeq\mathscr A_2(F'-,-),
    \]
    using the adjunction $F_1\dashv R_1$ and the full faithfulness of $F_2$.
    Finally, we can define $R'\coloneq R_1F_2$ and get $R'R_2\simeq R_1$ by uniqueness of adjoints.
  \item {``\ref{it:R-1-factors}$\implies$\ref{it:R-1-inverts-epsilon-2}'':}
    Assume that there is $R'$ with $R'R_2\simeq R_1$.
    Since $F_2$ is fully faithful, the unit $\eta_2\colon \id \to R_2F_2$ is invertible,
    hence also $\eta_2R_2$.
    Since $\eta_2R_2$ is a section of $R_2\varepsilon_2$ by the triangle identity,
    we conclude that $R_2\varepsilon_2$ is invertible, hence also $R_1\varepsilon_2\simeq R'R_2\varepsilon_2$.
  \item {``\ref{it:R-1-inverts-epsilon-2}$\implies$\ref{it:F-1-maps-to-image-of-F-2}'':}
    Assume that each component of $R_1\varepsilon_2$ is an equivalence.
    Let $a_1$ be an object of $\mathscr A_1$. For each object $X$ of
    $\mathscr C$, the map
    \begin{equation}
      \label{eq:colocality-comparison-map-}
      (\varepsilon_2X\circ -)\colon
      \mathscr C(F_1a_1,F_2R_2X)\xrightarrow{\simeq}\mathscr C(F_1a_1,X)
    \end{equation}
    is an equivalence, because after applying the adjunction
    $F_1\dashv R_1$, it is identified with
    $(R_1\varepsilon_2X\circ -)$.
    Let $Z\coloneq F_1a_1$;
    we must show that it lies in the essential image of $F_2$.
    Let $s\colon Z\to F_2R_2Z$ be a preimage up to homotopy of $\id_Z$
    under the equivalence \eqref{eq:colocality-comparison-map-}
    (specifically with $X\coloneqq Z$),
    so that we have $\varepsilon_2Z\circ s\simeq\id_Z$.
    It remain to show that there is also a homotopy $s\circ\varepsilon_2Z\simeq\id_{F_2R_2Z}$,
    because then we have $s\colon Z\simeq F_2R_2Z$,
    so that $Z$ lies in the essential image of $F_2$.
    By construction, both sides are homotopic after postcomposing with $\varepsilon_2 Z$.
    Moreover, the composite
    \[
      \mathscr C(R_2Z,R_2Z) \xrightarrow{F_2}
      \mathscr{C}(F_2R_2Z,F_2R_2Z)
      \xrightarrow{(\varepsilon_2 Z\circ -)} \mathscr{C}(F_2R_2Z, Z)
    \]
    is the adjunction equivalence,
    and the left map is an equivalence by the full faithfulness of $F_2$.
    We conclude that the postcomposition map $(\varepsilon_2 Z\circ -)$
    is an equivalence, so that we may cancel it from the homotopy to conclude the proof.
  \end{itemize}

  The final statement of the lemma follows directly from the construction of
  $F'$ and $R'$ and the adjunction $F'\dashv R_1F_2$.
\end{proof}

We obtain an equivalent reformulation of \ref{A2} in presence of \ref{A1}.

\begin{lem}
\label{lem:A2-full-subcats}
Given a family of adjunctions
$\{\adjunction{J_n}{\mathscr{C}_n}{\mathscr{C}}{R_n}\}_{n\geq 0}$
satisfying \ref{A1}, Condition \ref{A2} can equivalently be replaced with the following condition:
\begin{enumerate}[label={(A2')}, ref={(A2')}]
\item\label{A2p}
For each $n>0$, we have an inclusion
\[
J_{n-1}(\mathscr C_{n-1}) \subseteq J_n(\mathscr C_n)
\]
as full subcategories of $\mathscr C$.
\end{enumerate}
\end{lem}

\begin{proof}
Apply \cref{lem:factor-adjunction-ff} to $F_1=J_{n-1}$ and $F_2=J_n$: condition \ref{A2p} is condition \ref{it:F-1-maps-to-image-of-F-2} of that lemma, while \ref{A2} is condition \ref{it:R-1-inverts-epsilon-2}.
\end{proof}

\begin{lem}
\label{lem:adjunction-to-fin-stages}
Given a preinductive tower
$\{\adjunction{J_n}{\mathscr{C}_n}{\mathscr{C}}{R_n}\}_{n\geq 0}$
and $n>0$, there is an adjunction
\[
\adjunction{R_nJ_{n-1}}{\mathscr C_{n-1}}{\mathscr C_n}{R_{n-1}J_n}.
\]
\end{lem}

\begin{proof}
Apply \cref{lem:factor-adjunction-ff} to $F_1=J_{n-1}$ and $F_2=J_n$. By \ref{A1}, $J_n$ is fully faithful, and \ref{A2} is condition \ref{it:R-1-inverts-epsilon-2} of that lemma. Its final statement gives the desired adjunction.
\end{proof}

\begin{lem}
\label{Induction}
Given a preinductive tower 
$\{\adjunction{J_n}{\mathscr{C}_n}{\mathscr{C}}{R_n}\}_{n\geq0}$,
for all $m,d\geq0$ and $X$ in $\mathscr{C}$, we have
\[
R_mJ_{m+d}R_{m+d}X\simeq R_mX.
\]
\end{lem}

\begin{proof}
By \ref{A1}, each $J_n$ is fully faithful, so we may identify $\mathscr C_n$ with the full subcategory $J_n(\mathscr C_n)$ of $\mathscr C$. By \ref{A2p}, we have inclusions
\[
J_{n-1}(\mathscr C_{n-1})\subseteq J_n(\mathscr C_n)
\]
as full subcategories of $\mathscr C$. Composing these, we obtain an inclusion
\[
J_m(\mathscr C_m)\subseteq J_{m+d}(\mathscr C_{m+d})
\]
as full subcategories of $\mathscr C$. Applying \cref{lem:factor-adjunction-ff} to $F_1=J_m$ and $F_2=J_{m+d}$ gives
\[
R_mJ_{m+d}R_{m+d}X\simeq R_mX,
\]
as desired.
\end{proof}

\begin{proof}[Proof of \cref{prop:preinductive-reflective-subcats}]
Note that the datum of a fully faithful functor
$J_n\colon\mathscr{C}_n\to \mathscr{C}$ (with fixed codomain $\mathscr{C}$) is the same datum as a full subcategory $\mathscr{C}'_n\subseteq \mathscr{C}$ by replacing $J_n$ with its essential image. (A more precise way of saying this is that the inclusion from the poset of full subcategories of $\mathscr{C}$ into the full subcategory of the slice $\infty\mathscr{C}\mathrm{at}_{/\mathscr{C}}$ spanned by the fully faithful functors is an equivalence.)
The functor $J_n$ having an adjoint $R_n$ corresponds to the property that this full subcategory is coreflective.
Finally, \cref{lem:A2-full-subcats} says that \ref{A2} corresponds to the property that these full subcategories are nested.

To prove the last statement, we have to see that each functor $R_{n-1}J_n$ does indeed correspond to $R_{n-1}^n$. Note that the functor $R_nJ_{n-1}$ corresponds to the inclusion $\mathscr{C}_{n-1}\hookrightarrow \mathscr{C}_{n}$ because of the commutative diagram
\[\begin{tikzcd}
      \mathscr{C}_n
      &
      \mathscr{C}_{n-1}
      \ar[l,"R_nJ_{n-1}"]
      \ar[d,"J_{n-1}","\simeq"']
      \\
      J_n(\mathscr{C}_n)
      \ar[u,"R_n","\simeq"']
      \ar[r,hookleftarrow]
      \ar[d,hookrightarrow]
      &
      J_{n-1}(\mathscr{C}_{n-1})
      \\
      \mathscr{C}
      \ar[ur,hookleftarrow]
    \end{tikzcd}
\]
in $\infty\mathscr{C}\mathrm{at}$; passing to right adjoints yields the claimed correspondence, because the right adjoint of $R_nJ_{n-1}$ is $R_{n-1}J_{n}$ by \cref{lem:adjunction-to-fin-stages}.
\end{proof}

\begin{const}
\label{const:functor-R-vec}
Given a preinductive tower $\{\adjunction{J_n}{\mathscr{C}_n}{\mathscr{C}}{R_n}\}_{n\geq0}$,
the translation from \cref{prop:preinductive-reflective-subcats} gives rise to a commutative diagram
\begin{equation}
\label{eq:all-C-to-infty}
\begin{tikzcd}
      \mathscr{C}_0
      \ar[r,hookrightarrow]
      \ar[drrr,hookrightarrow]
      &
      \mathscr{C}_1
      \ar[r,hookrightarrow]
      \ar[drr,hookrightarrow]
      &
      \cdots
      \ar[r,hookrightarrow]
      \ar[dr, bend left =10,phantom,"\cdots"]
      &
      \mathscr{C}_{n-1}
      \ar[r,hookrightarrow]
      \ar[d,hookrightarrow]
      &
      \mathscr{C}_n
      \ar[r,hookrightarrow]
      \ar[dl,hookrightarrow]
      &
      \mathscr{C}_{n+1}
      \ar[r,hookrightarrow]
      \ar[dll,hookrightarrow]
      &
      \cdots
      \ar[dlll, bend left=15,phantom,"\cdots" near start]
      \\
      &&& \mathscr{C}
    \end{tikzcd}
\end{equation}
This diagram can be regarded as a functor $\mathscr{C}_\bullet\colon\mathbb{N}\cup\{\infty\}\to \infty\mathscr{C}\mathrm{at}$, where we view $\mathbb{N}\cup\{\infty\}$ as a poset with the usual order, hence a category, hence an $\infty$-category.
\end{const}

\begin{notn}
  Given a preinductive tower $\{\adjunction{J_n}{\mathscr{C}_n}{\mathscr{C}}{R_n}\}_{n\geq0}$,
  we denote by $p_{\mathscr{C}}\colon \int_{\mathbb N\cup\{\infty\}}\mathscr C_\bullet\to \mathbb N\cup\{\infty\}$
the Grothendieck construction of the functor $\mathscr{C}_\bullet\colon\mathbb{N}\cup\{\infty\}\to \infty\mathscr{C}\mathrm{at}$
as in \cite[\textsection3.2]{htt}.
Given a full sub-poset $B\subseteq \mathbb{N}\cup\{\infty\}$, we denote by
\[
  \mathrm{Sec}_B(p_{\mathscr{C}})\coloneqq
  \Fun_{\mathbb{N}\cup\{\infty\}}\left(B,\int_{\mathbb{N}\cup\{\infty\}}\mathscr C_\bullet\right)
  \simeq
  \Fun_B\left(B,\int_{B}\mathscr{C}_\bullet\right)
\]
the $\infty$-category of sections
\[\begin{tikzcd}
      &
      \int_{\mathbb N\cup\{\infty\}}
      \mathscr{C}_\bullet
      \ar[d,"{p_{\mathscr{C}}}"]
      \\
      B\ar[r,hookrightarrow]
      \ar[ur,dashed]
      &
      \mathbb{N}\cup\{\infty\}
    \end{tikzcd}
\]
of $p_{\mathscr{C}}$ on $B$.
For each further full sub-poset $A\subseteq B$,
we denote by
$\SecCart{B}{A}(p_{\mathscr{C}})\subseteq \mathrm{Sec}_B(p_{\mathscr C})$
the full subcategory of those sections
that are cartesian when restricted to $A$.
\end{notn}

\begin{lem}
\label{lem:emb-into-functors}
Fix $A\subseteq B\subseteq \mathbb N \cup \{\infty\}$.
The natural transformation $\mathscr{C}_\bullet \hookrightarrow \mathscr{C}$ yields a functor
\[\SecCart{B}{A}(p_{\mathscr{C}}) \hookrightarrow \mathscr{C}^B.\]
This functor is fully faithful and identifies the left side with the full subcategory of those $B$-indexed diagrams $\vec X$ for which
\begin{itemize}[leftmargin=*]
\item for each $n\in B$, the object $X_n$ lies in the full subcategory $\mathscr{C}_n \subseteq \mathscr{C}$ (with $\mathscr{C}_\infty\coloneq \mathscr{C}$) and
\item for each $m,n\in A$ with $m\leq n$,
  the induced map $X_m\to R_m^nX_n$ is an equivalence
  (where $R_m^n\colon \mathscr C_n\to \mathscr C_m$
  is a right adjoint to the inclusion $\mathscr C_m\subseteq \mathscr C_n$).
\end{itemize}
\end{lem}
 
\begin{proof}
After unstraightening, the natural transformation $\mathscr{C}_\bullet \hookrightarrow \mathscr{C}$ (of diagrams $\mathbb N\cup \{\infty\}\to \infty\mathscr{C}\mathrm{at}$) corresponds to a functor
\[
    \begin{tikzcd}
      \int_{\mathbb{N}\cup\{\infty\}}\mathscr{C}_\bullet
      \ar[r]
      \ar[dr,"{p_{\mathscr{C}}}"]
      &
      \int_{\mathbb{N}\cup\{\infty\}}\mathscr{C}
      \ar[r,"\simeq"]
      &
      (\mathbb{N}\cup\{\infty\})\times \mathscr{C}
      \ar[dl]
      \\
      &
      \mathbb N\cup\{\infty\}
    \end{tikzcd}
\]
over $\mathbb N\cup\{\infty\}$.
The desired functor is just given by postcomposition of sections with this functor, noting that a $B$-section of the projection on the right side is just a functor $B\to \mathscr{C}$. It is fully faithful because each transition map in the diagram \eqref{eq:all-C-to-infty} is fully faithful.
\end{proof}

\begin{lem}
\label{lem:description-of-lim}
Let $\{\adjunction{J_n}{\mathscr{C}_n}{\mathscr{C}}{R_n}\}_{n\geq0}$ be a preinductive tower. There is an equivalence
\[\lim_{n\geq0}^R\mathscr{C}_n \simeq \SecCart{\mathbb{N}}{\mathbb N}(p_{\mathscr{C}})\]
that, under the fully faithful embedding of \cref{lem:emb-into-functors},
identifies the right limit with the full subcategory of \(\mathscr C^{\mathbb{N}}\)
spanned by sequences of the form
\[\vec{X} : J_0X_0\to J_1X_1 \to \cdots \to J_nX_n \to \cdots\]
in $\mathscr{C}$ where each $X_n\in \mathscr{C}_n$ and such that for each $n\in \mathbb{N}$ the induced map
\[X_n \xrightarrow{\simeq} R_nJ_{n+1}X_{n+1}\]
is an equivalence.
\end{lem}

\begin{proof}
Note that the cocartesian fibration $p_{\mathscr{C}}$ is also cartesian because each inclusion $\mathscr{C}_m\hookrightarrow \mathscr{C}_n$ (including the case $n=\infty$) has a right adjoint, by the dual of \cite[Cor.~5.2.2.5]{htt}. The claimed equivalence 
\[\lim_{n\geq0}^R\mathscr{C}_n \simeq \SecCart{\mathbb{N}}{\mathbb N}(p_{\mathscr{C}})\]
now follows from \cite[Cor.~3.3.3.2]{htt}.
\end{proof}

\begin{prop}
\label{FormulaGamma}
Every preinductive tower
$\{\adjunction{J_n}{\mathscr{C}_n}{\mathscr{C}}{R_n}\}_{n\geq0}$
induces an adjunction
\[
\adjunction{\gamma}{\lim_{n\geq0}^{R}\mathscr{C}_n}{\mathscr{C}}{\vec R}.
\]
Using the explicit description of \cref{lem:description-of-lim} for objects of the limit,
the two adjoints are described by the following formulas:
\[
\gamma(\vec X)\coloneqq\colim_{\mathbb{N}}\vec X
\]
and
\[\vec R Y\simeq (J_0R_0 Y\to J_1R_1Y \to \cdots \to J_nR_nY\to J_{n+1}R_{n+1}Y\to\cdots)
\]
with the transition maps adjoint to the equivalences
$R_nY \simeq R_n J_{n+1}R_{n+1}Y$
obtained from Condition \ref{A2} by passing to inverse.
\end{prop}

\begin{proof}
By \cite[\textsection4.3]{htt} the map $j\colon \mathbb N\hookrightarrow\mathbb N\cup\{\infty\}$ induces an adjunction
\[\adjunction{j_!}{\mathrm{Sec}_{\mathbb N}(p_{\mathscr{C}})}
  {\mathrm{Sec}_{\mathbb N\cup\{\infty\}}(p_{\mathscr{C}})}{j^*},
\]
where the left adjoint $j_!$ is given by relative left Kan extension and is fully faithful. Under the respective embeddings of \cref{lem:emb-into-functors},
the objects of these $\infty$-categories are described, respectively, as sequences
\[\vec{X} =
(J_0X_0\to J_1X_1 \to \cdots \to J_nX_n\to J_{n+1}X_{n+1}\to\cdots),\]
and
\[\vec{Y}=(J_0Y_0\to J_1Y_1 \to \cdots \to J_nY_n\to J_{n+1}Y_{n+1}\to\cdots\to Y),\]
with $X_n, Y_n$ objects of $\mathscr C_n$ and $Y$ an object of $\mathscr C$. The functors $j_!$ and $j^*$ are described by
\[j_!(\vec X)\simeq(J_0X_0\to J_1X_1 \to \cdots \to J_nX_n\to J_{n+1}X_{n+1}\to\cdots\to \colim_{n\geq0}J_nX_n),\]
(extending by the colimit at position $\infty$)
and
\[j^*\vec Y\simeq
(J_0Y_0\to J_1Y_1 \to \cdots \to J_nY_n\to J_{n+1}Y_{n+1}\to\cdots)\]
(restricting the section to $\mathbb{N}\subseteq \mathbb{N}\cup\{\infty\}$), respectively.
Moreover it is clear that the adjunction $(j_!,j^*)$ restricts and corestricts to an adjunction
\begin{equation}
\label{eq:adjunction-seccartN+1}
\adjunction{j_!}{\SecCart{\mathbb N}{\mathbb N}(p_{\mathscr{C}})}{\SecCart{\mathbb N\cup\{\infty\}}{\mathbb N}(p_{\mathscr{C}})}{j^*}.
\end{equation}

By \cite[\textsection4.3]{htt},
the map $i\colon \{\infty\}\hookrightarrow\mathbb N\cup\{\infty\}$ induces an adjunction
\[
\adjunction{i^*}{\mathrm{Sec}_{\mathbb N\cup\{\infty\}}(p_{\mathscr{C}})}
  {\mathscr C_\infty}{i_*},
\]
with the right adjoint $i_*$ fully faithful. Here, the right adjoint functor $i_*$ acts on objects as
\[
i_*(X)\cong(J_0R_0X\to J_1R_1X \to \cdots \to J_nR_nX\to J_{n+1}R_{n+1}X\to\cdots\to X),
\]
and the left adjoint functor $i^*$ acts on objects as
\[
i^*\vec Y\cong Y.
\]
By Condition \ref{A2}, the section $i_*(X)$ takes values in cartesian morphisms when restricted to $\mathbb N$.
In particular, the adjunction restricts to
\begin{equation}
\label{eq:adjunction-seccart-C}
\adjunction{i^*}{\SecCart{\mathbb N\cup\{\infty\}}{\mathbb N}(p_{\mathscr{C}})}{\mathscr C_\infty}{i_*}.
\end{equation}

Composing the two adjunctions
\eqref{eq:adjunction-seccartN+1} and \eqref{eq:adjunction-seccart-C}
with the equivalence of \cref{lem:description-of-lim} we obtain the desired adjunction
\[
  \adjunction{\gamma\coloneqq i^*j_!}{\lim_{n\geq 0}^R\mathscr{C}_n}{\mathscr C_\infty}{j^*i_*\eqqcolon\vec R},\]
satisfying the desired formulas.
\end{proof}

\begin{defn}
\label{defn-axioms-inductive}
Let $\{\adjunction{J_n}{\mathscr{C}_n}{\mathscr{C}}{R_n}\}_{n\geq0}$ be a preinductive tower. We introduce the following conditions:
\begin{enumerate}[label=(A\arabic*), ref=(A\arabic*), start=3]
\item \label{A3}
For all $m\geq 0$ and all $X$ in $\mathscr{C}$, the canonical map
\[\colim_{n\geq 0}R_mJ_nR_nX\to R_m\colim_{n\geq 0}J_nR_nX\]
is an equivalence.
\item \label{A4}
The functors $R_n$ for $n\geq 0$ are jointly conservative; that is, any $f\colon X\to Y$ in $\mathscr{C}$ which induces an equivalence
\[R_nf\colon R_nX\xrightarrow{\simeq} R_nY\]
for all $n\geq 0$ is itself an equivalence $f\colon X\xrightarrow{\simeq} Y$.
\item \label{A5}
Any $\vec f\colon\vec X\to\vec Y$ in $\lim_{n\geq 0}^{R}\mathscr{C}_n$ which induces an equivalence
\[\colim_{n\geq 0}J_nf_n\colon\colim_{n\geq 0}J_nX_n\xrightarrow{\simeq}\colim_{n\geq 0}J_nY_n\]
in $\mathscr{C}$ is itself an equivalence $\vec f\colon\vec X\xrightarrow{\simeq}\vec Y$.
\end{enumerate}
We also consider the following variants:
\begin{enumerate}[label={(A3')}, ref={(A3')}]
\item\label{A3p} For each $m\geq 0$ and each $X$ in $\mathscr C$,
  the canonical map
\[R_m X \xrightarrow[\simeq]{\eta_mR_m X} R_mJ_mR_m X\to R_m\colim_{n\geq 0}J_nR_nX,\]
induced by the colimit structure map
$J_mR_mX\to\colim_{n\geq 0}J_nR_nX$, is an equivalence.
\end{enumerate}
\begin{enumerate}[label=(B\arabic*), ref=(B\arabic*), start=3]
\item \label{B3} For each $X$ in $\mathscr C$, the canonical map
\[\colim_{n\geq 0}J_nR_nX \to X,\]
induced by the compatible system
of counit maps $\varepsilon_{n}X\colon J_nR_nX\to X$, is an equivalence.
\item\label{B4} For each $\vec{X}$ in $\lim^R_{n\geq 0}\mathscr{C}_n$ and each $m\in \mathbb{N}$, the canonical map
\[
X_m \to R_m\colim_{n\geq 0}J_n X_n
\]
(adjoint to the structure map of $J_m X_m$ into the colimit)
is an equivalence.
\end{enumerate}
\end{defn}

\begin{lem}
\label{ReformulationA3}
For a preinductive tower,
conditions \ref{A3} and \ref{A3p} are equivalent.
\end{lem}

\begin{proof}
By \cref{Induction} there is a canonical equivalence
\[
\colim_{n\geq 0}R_mJ_nR_nX \simeq R_mX.
\]
Under this equivalence, the canonical map
\[
\colim_{n\geq 0}R_mJ_nR_nX\to R_m\colim_{n\geq 0}J_nR_nX
\]
from Condition \ref{A3} is identified with the map
in Condition \ref{A3p}. Hence \ref{A3} and \ref{A3p} are equivalent.
\end{proof}

We can tie the conditions of \cref{defn-axioms-inductive} directly to properties
of the adjunction of \cref{FormulaGamma}:

\begin{thm}
\label{ReformulationB3}
Consider a preinductive tower
$\{\adjunction{J_n}{\mathscr{C}_n}{\mathscr{C}}{R_n}\}_{n\geq0}$
and the corresponding adjunction
\[
\adjunction{\gamma}{\lim_{n\geq0}^{R}\mathscr{C}_n}{\mathscr{C}}{\vec R}.
\]
\begin{enumerate}[leftmargin=*]
\item Condition \ref{A4} says precisely that the right adjoint $\vec{R}$ is conservative.
\item Condition \ref{A5} says precisely that the left adjoint $\gamma$ is conservative.
\item Condition \ref{B3} says exactly that the right adjoint $\vec{R}$ is fully faithful. Moreover, Condition \ref{B3} is equivalent to the conjunction (\ref{A3} and \ref{A4}).
\item Condition \ref{B4} says exactly that the left adjoint $\gamma$ is fully faithful.
\end{enumerate}
\end{thm}

\begin{proof}
\begin{enumerate}[leftmargin=*]
\item
On its face, Condition \ref{A4} says that the composite
\[\mathscr{C}\xrightarrow{\vec{R}} \lim_{n\geq 0}^R\mathscr{C}_n \to \prod_{n\geq 0}\mathscr{C}_n\]
is conservative. But this is the same as saying that $\vec{R}$ is conservative, because the map from the limit to the product is conservative.
\item
  There is nothing to prove, because the induced map in \ref{A5} is precisely $\gamma(\vec {f})$.
\item On its face Condition \ref{B3} says that the counit of the adjunction is an equivalence;
  this is equivalent to saying that the right adjoint is fully faithful
  which in turn implies that the right adjoint is conservative, i.e., Condition \ref{A4}.

We conclude the proof by showing that Condition \ref{A3} and \ref{B3} are equivalent in the presence of Condition \ref{A4}. By \cref{ReformulationA3} we may work with Condition \ref{A3p} instead of Condition \ref{A3}.
    
Assume Condition \ref{A4}. Let $c_X\colon \colim_{n\geq 0}J_nR_nX\to X$
denote the map in Condition \ref{B3}.
For each $m\geq 0$, the composite
\[
  \begin{tikzcd}[column sep=large]
    R_mX\ar[r,"\eta_m R_m X"]
    &
    R_mJ_mR_mX\ar[r]
    \ar[rr,bend left=10,"R_m\varepsilon_m X"]
    &R_m\left(\colim_{n\geq 0}J_nR_nX\right)
    \ar[r,"R_mc_X"']
    &
    R_mX
  \end{tikzcd}
\]
is the identity by the triangle identity for the adjunction $J_m\dashv R_m$.
Thus $R_mc_X$ is a retraction of the composite map in \ref{A3p}, so one is an equivalence if and only if the other is. By \ref{A4}, $c_X$ is an equivalence if and only if $R_mc_X$ is an equivalence for every $m\geq 0$. Therefore \ref{A3p} and \ref{B3} are equivalent in the presence of Condition \ref{A4}.
\item On its face Condition \ref{B4} says that the unit $\vec{X}\to \vec{R}\gamma(\vec{X})$ becomes an equivalence after evaluating at every $m\in \mathbb{N}$.
  This suffices for the unit to be an equivalence (i.e., the left adjoint $\gamma$ being fully faithful),
  again because the projections from the limit to its stages are jointly conservative.\qedhere
\end{enumerate}
\end{proof}

As an immediate corollary, we get the following recognition principle for right towers:

\begin{cor}
\label{MainTheorem}
Given a preinductive tower
$\{\adjunction{J_n}{\mathscr{C}_n}{\mathscr{C}}{R_n}\}_{n\geq0}$
the following are equivalent:
\begin{itemize}[leftmargin=*]
\item The adjunction $\adjunction{\gamma}{\lim_{n\geq0}^{R}\mathscr{C}_n}{\mathscr{C}}{\vec R}$ from \cref{FormulaGamma}
is an equivalence.
\item Conditions \ref{A3}, \ref{A4}, and \ref{A5} hold.
\item Conditions \ref{B3} and \ref{B4} hold.
\end{itemize}
\end{cor}

\begin{proof}
This follows immediately from the
fact that an adjunction is an equivalence
if and only if one adjoint is fully faithful and the other is conservative.
\end{proof}

\begin{defn}
An \emph{inductive tower} $\{\adjunction{J_n}{\mathscr{C}_n}{\mathscr{C}}{R_n}\}_{n\geq0}$ is a preinductive tower that satisfies the equivalent conditions of \cref{MainTheorem}.
\end{defn}

\section{Limit of strict \texorpdfstring{$n$}{n}-categories}

\label{Ncategories}
We recall that $\omega\cat$ denotes the category\footnote{In this section we are really considering the ordinary categories $\omega\cat$ and $n\cat$ themselves, \emph{not} the $\infty$-categories they present through the model structures from \cite{LMW}.}
of strict $\omega$-categories and,
for all $n\geq0$,
let $n\cat$ denote the category of strict $n$-categories,
as in e.g.~\cite[\textsection3]{LMW} or  \cite[\textsection1.2]{AraMaltsiniotisJoin}.
Under the canonical fully faithful embedding $I_n\colon n\cat\hookrightarrow\omega\cat$,
we identify each strict $n$-category with the corresponding strict $\omega$-category
(with the same arrows up to dimension $n$ and only identities above it.)

\subsection{Right limit: Core filtration of strict \texorpdfstring{$\omega$}{omega}-categories}
Recall 
that, for all $n\geq0$, there is an adjunction
\[
\adjunction{I_n}{n\cat}{\omega\cat}{\mathrm{core}_n},
\]
where $\mathrm{core}_n(X)$ is the $n$-th \emph{core}, a.k.a.\ \emph{naive truncation}, of $X$, which satisfies
\[
b\Mor(\mathrm{core}_nX)=b\Mor(X)\quad\text{for $0\leq b\leq n$}.
\]

\begin{prop}
\label{Verification-omega-cats-right}
The family of adjunctions
$\{\adjunction{I_n}{n\cat}{\omega\cat}{\mathrm{core}_n}\}_{n\geq0}$
is an inductive tower for $\omega\cat$.
\end{prop}

\begin{proof}
  For each $n\geq 0$, the image $I_n(n\cat)$ is the full subcategory of $\omega\cat$
  consisting of those strict $\omega$-categories
  that have only identity morphisms in dimensions $>n$.
  These full subcategories are clearly nested and have coreflectors $I_n\core_n$,
  hence \ref{A1} and \ref{A2}
  follow from \cref{prop:preinductive-reflective-subcats}.

  We verify \ref{A3} and \ref{A4}
  (which are jointly equivalent to \ref{B3}) and \ref{B4}:

For \ref{A3}, the claim follows from the fact that each functor $\mathrm{core}_m$ preserves all colimits, because it has a further right adjoint by \cite[\textsection1.2]{AraMaltsiniotisJoin}.

For \ref{A4}, we observe that, if $\mathrm{core}_nf$ is an isomorphism for all $n\geq0$, then $n\Mor(\mathrm{core}_nf)=n\Mor f$ is an isomorphism for all $n$, so $f$ is an isomorphism.

For \ref{B4}, let
\[
X_0\to X_1\to\cdots\to X_n\to\cdots\to\colim_n X_n
\]
be a colimit sequence that induces isomorphisms
\[
X_n\cong \mathrm{core}_nX_n\xrightarrow{\cong}\mathrm{core}_nX_{n+1}.
\]
Then, for any fixed $m$, applying $\core_m$ yields a colimit sequence (again, because $\core_m$ preserves colimits) which stabilizes above level $m$, hence induces an isomorphism
\[
X_m\xrightarrow{\cong}\core_m(\colim_n X_n),
\]
concluding the proof.
\end{proof}

\begin{cor}
\label{LimitIndStrict}
There is an equivalence of categories $\gamma\colon\lim^{R}_{n\geq0} n\cat \simeq \omega\cat\noloc\vec{\mathrm{core}}$.
\end{cor}

\begin{proof}
  By definition of inductive tower.
\end{proof}

\subsection{Left limit: Truncation tower of strict \texorpdfstring{$\omega$}{omega}-categories}
We recall that, for all $n\geq0$, there is an adjunction
\[
\adjunction{\mathrm{tr}_n}{\omega\cat}{n\cat}{I_n},
\]
where $\mathrm{tr}_n(X)$ is the $n$-th \emph{intelligent truncation} of $X$,
which satisfies
\[
b\Mor(\mathrm{tr}_nX)=
\begin{cases}
b\Mor(X)&b<n\\
n\Mor(X)/(n+1)\Mor(X)&b=n.
\end{cases}
\]
Here, we denote by $n\Mor(X)/(n+1)\Mor(X)$ the quotient of $n\Mor(X)$
modulo the equivalence relation generated by the condition $f\simeq g$
whenever there exists $\alpha\colon f\to g$ in $(n+1)\Mor(X)$.
See \cite[\textsection1.2]{AraMaltsiniotisJoin} for details.

We use the dual of the framework in \Cref{Framework}.

\begin{prop}
\label{Verification-omega-cats-left}
The family of adjunctions
$\{\adjunction{\mathrm{tr}_n}{\omega\cat}{n\cat}{I_n}\}_{n\geq0}$
is a dual inductive tower for $\omega\cat$.
\end{prop}

\begin{proof}
For the duals of \ref{A1} and \ref{A2}, the family is a dual preinductive tower because $\{n\cat\}_n$ is a system of nested full reflective subcategories of $\omega\cat$, by \cref{prop:preinductive-reflective-subcats}.

For the dual of \ref{B3}, we have to verify that for each strict $\omega$-category $X$, the sequence
\[
X\to\cdots\to\mathrm{tr}_nX\to\cdots\to\mathrm{tr}_1X\to\mathrm{tr}_0X
\]
is a limit sequence. Since limits of strict $\omega$-categories are computed levelwise on arrow sets, it suffices to observe that for each $n>a\geq0$ one has
\[
a\Mor(X)\xrightarrow{\cong}a\Mor(\tr_nX),
\]
so that the limit stabilizes levelwise and yields the desired value.

For the dual of \ref{B4}, consider a sequence
\[
\cdots\to X_n\to\cdots\to X_1\to X_0
\]
of strict $\omega$-categories which induces isomorphisms
\[
\tr_n(X_{n+1})\xrightarrow{\cong}X_n.
\]
We have to show that the induced map
\[
\tr_m(\lim_nX_n)\to X_m
\]
is an isomorphism for all $m$. Note that, for fixed $a\leq m+1$, the sequence
\[
\cdots\to a\Mor(X_n)\to\cdots\to a\Mor(X_1)\to a\Mor(X_0)
\]
stabilizes for $n\geq m+2$, hence we have isomorphisms
\[
a\Mor(\lim_nX_n)\xrightarrow{\cong}a\Mor(X_{m+2}).
\]
Since the $m$-truncation only depends on the arrow sets up to dimension $m+1$, this means that we have
\[
\tr_m(\lim_nX_n)\xrightarrow{\cong}\tr_m(X_{m+2})\xrightarrow{\cong}X_m,
\]
as desired.
\end{proof}

\begin{cor}
\label{LimitCoIndStrict}
There is an equivalence of categories
$\vec{\mathrm{tr}}\colon\omega\cat\simeq\lim^{L}_{n\geq0} n\cat\noloc\lambda$,
with $\lambda(\vec X)=\lim_n I_nX_n$.
\end{cor}
\begin{proof}
  By definition of (dual) inductive tower.
\end{proof}

\begin{rmk}
One could have also directly constructed an equivalence
\[
\lim^R_nn\cat \simeq \lim^L_n n\cat
\]
between the left and the right limit.
The only thing one needs to use is that the $(n+1)$-truncation does not change the $n$-core, and that the $n$-truncation only depends on the $(n+1)$-core, so that one has
\[
\core_{n}\tr_{n+1}= \core_n
\quad\text{and}\quad
\tr_{n}\core_{n+1}=\tr_n.
\]
Indeed, using these identities one can compute the limit of the alternating sequence 
\[
0\cat \xleftarrow{\core_0} 1\cat
\xleftarrow{\tr_1}\cdots\leftarrow 2k\cat\xleftarrow{\core_{2k}}(2k+1)\cat\xleftarrow{\tr_{2k+1}}\cdots
\]
in two ways, using that both the even and odd numbers are cofinal in $(\mathbb{N},\leq)$: The even entries converge to $\lim^R_{k}2k\cat \simeq \lim^R_nn\cat$ while the odd entries converge to $\lim^L_{k}(2k+1)\cat \simeq \lim^L_nn\cat$.

Amusingly, this argument also shows that the limit $\lim^q_{n} n\cat$ is independent of the choice of transition functors $q_n\colon (n+1)\cat \to n\cat$ even if each $q_n\in\{\core_n,\tr_n\}$ is chosen independently.
Indeed there must exist a cofinal subset $A\subseteq \mathbb{N}$
such that either
$q_n=\core_n$ for all $n\in A$
or
$q_n=\tr_n$ for all $n\in A$.
Restricting to $A$ one then sees that $\lim^q_n n\cat$
is equivalent to either the right or left limit, respectively.
\end{rmk}

\section{Limit of \texorpdfstring{$(\infty,n)$}{(infinity,n)}-categories}
\label{ooNcategories}

We recall from \cite{VeritySlides,EmilyNotes,OR,RiehlVerityBook} that a
\emph{saturated complicial set}
is a marked simplicial set satisfying the right lifting property with respect
to a certain set of \emph{elementary anodyne extensions}.
These are recalled in \cite[Definition~1.19]{OR} and are:
\begin{enumerate*}
    \item the \emph{complicial horn extensions},
    \item the \emph{thinness extensions} and
    \stepcounter{enumi}
    \item the \emph{saturation extensions}.
\end{enumerate*} 
Additionally adding (3) the \emph{triviality extensions} in dimension above $n$
defines the \emph{$n$-complicial sets}.
There are model structures $m\sset_{\mathrm{cmp,sat}}$ and $m\sset_{(\infty,n)}$
whose fibrant objects are the saturated complicial sets
and the $n$-complicial sets, respectively;
in both cases, the cofibrations are the monomorphisms of underlying simplicial sets.
We denote by $\mathscr Comp$ and $\mathscr Comp_n$
the $\infty$-categories presented by these model structures, respectively.
\begin{rmk}
\label{defn:Bousfield-Comp}
The Quillen adjunction
\[
\adjunction{\id}{m\sset_{\mathrm{cmp,sat}}}{m\sset_{(\infty,n)}}{\id}
\]
presents an adjunction
\[\adjunction{L_n\coloneqq \mathbb{L}\id}{\mathscr Comp}{\mathscr Comp_n}{\mathbb{R}\id\eqqcolon I_n}\]
of $\infty$-categories.
It is a reflective localization because the strict counit
$\id\circ \id \xrightarrow{=} \id$
already computes the derived counit $L_nI_n\to \id$
(because every object is cofibrant), which is therefore an equivalence.
Thus the right adjoint $I_n\colon \mathscr Comp_n\hookrightarrow\mathscr Comp$
is fully faithful and its image consists of those objects of
$\mathscr Comp$ that admit an $n$-complicial representative in $m\sset$.
\end{rmk}

We collect here a couple of properties of complicial sets that we will need later.

\begin{rmk}
\label{SpPreservesFilteredColimits}
As discussed in \cite[Remark~3.3]{ORW2}, for every $n\geq0$, the functor $\mathrm{sp}_n\colon m\sset\to m\sset$
preserves filtered colimits.
\end{rmk}

\begin{lem}
\label{LemmaReconstruction}
Let
\[
X_0\subseteq X_1\subseteq\dots\subseteq X_n\subseteq X_{n+1}\subseteq\dots
\]
be a sequence of marked simplicial sets such that $X_n$ is fibrant in $m\sset_{(\infty,n)}$ for every $n\geq0$. Then $\colim_{n\geq0}X_n$ is fibrant in $m\sset_{\mathrm{cmp,sat}}$.
\end{lem}

\begin{proof}
  Let $I\to J$ be an elementary anodyne extension of $m\sset_{\mathrm{cmp,sat}}$
  (complicial horn, thinness, or saturation).
  In all cases,
  $I$ has finitely many nondegenerate simplices.
  Consider a map
    $I\to\colim_{n\geq0}X_n$.
  By finiteness, this map factors through
    $I\to X_N$
  for some $N\geq0$. Since $X_N$ is fibrant in $m\sset_{(\infty,N)}$, it is in particular fibrant in $m\sset_{\mathrm{cmp,sat}}$. Hence there exists a lift
    $J\to X_N$.
  Postcomposing with $X_N\to\colim_{n\geq0}X_n$ gives a lift
    $J\to\colim_{n\geq0}X_n$.
  Thus $\colim_{n\geq0}X_n$ has the right lifting property with respect to all elementary anodyne extensions of $m\sset_{\mathrm{cmp,sat}}$, and is therefore fibrant.
\end{proof}

\subsection{Right limit: Core tower of a saturated complicial set}

We construct a preinductive tower by introducing the further right adjoint of $I_n$.

\begin{prop}
\label{Adjunction}
For all $n\geq0$, there is an adjunction
\[
\adjunction{I_n}{\mathscr Comp_n}{\mathscr Comp}{\mathrm{core}_n},
\]
where $\mathrm{core}_n(X)$ is the $n$-th \emph{core},
which essentially forgets the non-invertible $(n+1)$-morphisms.
\end{prop}

\begin{proof}
There is the Quillen pair
\[
\adjunction{\mathrm{th}_n}{m\sset_{(\infty,n)}}{m\sset_{\mathrm{cmp,sat}}}{\mathrm{sp}_n}
\]
from \cite[Prop.~4.5]{ORW2}.
We claim that $\mathbb{L}\mathrm{th}_n\simeq I_n$,
so that the claimed right adjoint is given by $\mathbb{R}\mathrm{sp}_n$.
Since $I_n$ is the right derived functor of the identity functor,
it suffices to observe that we have literal equality $X= \mathrm{th}_n(X)$
for any $n$-complicial set $X$.
\end{proof}

\begin{thm}
\label{PreinductiveComplicial}
The family of adjunctions
$\{\adjunction{I_n}{\mathscr Comp_n}{\mathscr Comp}{\mathrm{core}_n}\}_{n\geq0}$
is a preinductive tower for $\mathscr Comp$.
\end{thm}

\begin{proof}
Each $I_n\colon\mathscr Comp_n\to\mathscr Comp$ is fully faithful (by \cref{defn:Bousfield-Comp}) and has the right adjoint $\mathrm{core}_n$ by \cref{Adjunction}. Moreover, the essential images are nested,
\[
I_n(\mathscr Comp_n)\subseteq I_{n+1}(\mathscr Comp_{n+1}),
\]
because every $n$-complicial set is also $(n+1)$-complicial.
Thus $\{\mathscr Comp_n\}_{n\geq0}$ is a nested system of full coreflective subcategories of $\mathscr Comp$, and the claim follows from \cref{prop:preinductive-reflective-subcats}.
\end{proof}

We will show in \cref{InductiveComplicial}
that this is indeed an inductive tower,
so that $\mathscr Comp$ is equivalent to the right limit of the $\mathscr Comp_n$.
For this computation we will need two basic strictification lemmas that we prove next.

Given a model category $\cM$, we denote by $[\cM]$ its associated $\infty$-category,
and by $[-]\colon \cM\to [\cM]$ the localization functor.

\begin{lem}
\label{RepresentingCofibration}
Let $\cM$ be a model category.
Given a map $g\colon X\to X'$ in the associated $\infty$-category $[\cM]$,
a cofibrant object $Y$ in $\cM$,
and an equivalence $h\colon [Y]\simeq X$, there exist a fibrant object $Y'$ in $\cM$, a cofibration $f\colon Y\hookrightarrow Y'$ in $\cM$, and a commutative square in $[\cM]$
\[
\begin{tikzcd}
{[Y]}\arrow[r,"{[f]}"]\arrow[d,"\simeq"swap,"h"]&{[Y']}\arrow[d,"\simeq"]\\
X\arrow[r,"g"]&X'.
\end{tikzcd}
\]
\end{lem}

This is a routine argument but we include it for completeness. 

\begin{proof}
There exists a fibrant object $Y''$ in $\cM$ together with an equivalence
$[Y'']\simeq X'$.
Since $Y$ is cofibrant in $\cM$ and $Y''$ is fibrant in $\cM$, the composite
\[
[Y]\xrightarrow{h}X\xrightarrow{g}X'\simeq [Y'']
\]
is represented by a map $s\colon Y\to Y''$ in $\cM$. Factor $s$ as a cofibration followed by an acyclic fibration:
\[
s\colon Y \stackrel{f}{\hookrightarrow} Y'
\stackrel{\simeq}{\twoheadrightarrow} Y''.
\]
Since $Y''$ is fibrant and $Y'\to Y''$ is a fibration, $Y'$ is fibrant. Moreover, $Y'\to Y''$ is a weak equivalence, so it induces an equivalence
$[Y']\simeq [Y'']$
in $[\cM]$. We therefore obtain a commutative diagram
\[
\begin{tikzcd}
    {[Y]}\arrow[r,"{[f]}"]\arrow[d,"\simeq"swap,"h"]
    &{[Y']}\arrow[r,"\simeq"]\arrow[d,"\simeq"]
    &{[Y'']}\arrow[dl,"\simeq"]\\
    X\arrow[r,"g"]&X',
\end{tikzcd}
\]
which gives the desired square.
\end{proof}

\begin{lem}
\label{LemmaUnit}
Given $X_\bullet\in\lim_n^R\mathscr Comp_n$, there exists $Y_\bullet\in m\sset^{\mathbb N}$ such that
\begin{itemize}
\item $Y_n$ is $n$-complicial for every $n\geq0$;
\item $Y_n\hookrightarrow Y_{n+1}$ is a cofibration for every $n\geq0$;
\item $[Y_\bullet]\simeq I_\bullet X_\bullet$ in the common ambient\footnote{
    Here $I_\bullet X_\bullet$ denotes the image of $X_\bullet$
    under the fully faithful embedding $\lim_n^R\mathscr Comp_n \hookrightarrow \mathscr Comp^{\mathbb N}$ of \cref{lem:description-of-lim,lem:emb-into-functors}.}
$\infty$-category $\mathscr Comp^{\mathbb N}$.
\end{itemize}
Moreover the induced map
\[
\colim_n I_n X_n \xrightarrow{\simeq} \left[ \colim_n Y_n \right]
\]
is an equivalence in $\mathscr Comp$,
where the left colimit is computed in $\mathscr Comp$
and the right one is computed (strictly) in $m\sset_{\mathrm{cmp,sat}}$.
\end{lem}

\begin{proof}
  Inductively, we construct for every $n\geq0$ an object
  $Y_{0\leq\bullet\leq n}\in m\sset^{\mathbb N_{\leq n}}$ such that
  \begin{itemize}
  \item
    $Y_j$ is $j$-complicial for $0\leq j\leq n$ and
  \item
    $Y_j\hookrightarrow Y_{j+1}$ is a cofibration for $0\leq j<n$,
  \item
    together with
    an equivalence $h_\bullet \colon [Y_{0\leq\bullet\leq n}]\simeq (I_\bullet X_\bullet)_{0\leq\bullet\leq n}$ in the common ambient $\infty$-category $\mathscr Comp^{\mathbb N_{\leq n}}$.
  \end{itemize}
For the base case, let $Y_0$ be a $0$-complicial set representing $X_0$ in
$\mathscr Comp_0$.
For the inductive step, suppose that
$Y_{0\leq \bullet \leq n}$ and $h_{0\leq \bullet \leq n}$
have already been constructed.
Denote by $I^{n+1}_n\colon \mathscr Comp_n\to\mathscr Comp_{n+1}$
the fully faithful inclusion.
Since $Y_n$ is also $(n+1)$-complicial it represents $I_n^{n+1}(X_n)$ in
$\mathscr Comp_{n+1}$.
Thus we have the equivalence
\[
h_n\colon [Y_n]_{n+1}\simeq I_n^{n+1}([Y_n]_n)\xrightarrow{\simeq} I_n^{n+1}(X_n),
\]
where for all $m\in \mathbb{N}$
we denote by $[-]_m\colon m\sset \to \mathscr Comp_m$ the localization functor to remove ambiguity.
Applying \cref{RepresentingCofibration} in the model category $m\sset_{(\infty,n+1)}$ to the cofibrant marked simplicial set $Y_n$, the equivalence $h_n$, and the structure map
\[
  I_n^{n+1}(X_n)\to X_{n+1},
\]
we obtain an $(n+1)$-complicial set $Y_{n+1}$, a cofibration $Y_n\hookrightarrow Y_{n+1}$, and a commutative square
\[
\begin{tikzcd}
{[Y_n]_{n+1}}\arrow[r,"{[f_n]}"]\arrow[d,"\simeq"swap,"h_n"]
    &{[Y_{n+1}]_{n+1}}\arrow[d,"\simeq"]\\
    I_n^{n+1}(X_n)\arrow[r]
    &X_{n+1}.
\end{tikzcd}
\]
This extends the equivalence of truncated diagrams from $\mathbb N_{\leq n}$ to $\mathbb N_{\leq n+1}$ and completes the induction. Passing to all $n$ gives the desired $Y_\bullet$.

For ``moreover'', we only have to show that the colimit $\colim_n Y_n$
is already a homotopy colimit in $m\sset_{\mathrm{cmp,sat}}$.
This follows from the fact that
$Y_0$ is (vacuously) cofibrant and
each $Y_n\hookrightarrow Y_{n+1}$ is a cofibration.
\end{proof}

We can now prove the theorem.

\begin{thm}
\label{InductiveComplicial}
The family of adjunctions
$\{\adjunction{I_n}{\mathscr Comp_n}{\mathscr Comp}{\mathrm{core}_n}\}_{n\geq0}$
is an inductive tower for $\mathscr Comp$.
\end{thm}

\begin{proof}
For \ref{B3}, let $X$ be a saturated complicial set. Since the inclusions
\[
\mathrm{sp}_0X\hookrightarrow\mathrm{sp}_1X\hookrightarrow\cdots\hookrightarrow\mathrm{sp}_nX\hookrightarrow\cdots
\]
are cofibrations and
\[
\colim_{n\geq0}\mathrm{sp}_nX\cong X
\]
as marked simplicial sets, the strict colimit represents the homotopy colimit in $m\sset_{\mathrm{cmp,sat}}$.
Hence the canonical map
\[
\colim_{n\geq0}I_n\mathrm{core}_n[X]\to [X]
\]
is an equivalence in $\mathscr Comp$, proving \ref{B3}.

For \ref{B4}, let $X_\bullet$ be an object of $\lim_{n\geq0}^{R}\mathscr Comp_n$.
By \cref{LemmaUnit}, we can choose $Y_\bullet\in m\sset^{\mathbb N}$
representing $I_\bullet X_\bullet$ such that each $Y_n$ is $n$-complicial and each structure map $Y_n\hookrightarrow Y_{n+1}$ is a cofibration.
Moreover, \cref{LemmaUnit} identifies the $m$-th component
$X_m\to \mathrm{core}_m(\colim_n I_nX_n)$
of the derived unit,
up to equivalence, with the strict map
\[
Y_m\to\mathrm{sp}_m\left(\colim_nY_n\right),
\]
using that $\colim_nY_n$ represents $\colim_n I_n X_n$
and is already fibrant in $m\sset_{\mathrm{cmp,sat}}$ (see \cref{LemmaReconstruction}).
It therefore suffices to show that this map is a weak equivalence in $m\sset_{(\infty,m)}$.
Since $\mathrm{sp}_m$ preserves filtered colimits (see \cref{SpPreservesFilteredColimits}),
we may exchange it with the colimit and instead consider the sequential composite
\[
Y_m =\mathrm{sp}_mY_m\to\mathrm{sp}_mY_{m+1}\to\ldots\to\colim_{n}\mathrm{sp}_mY_{n}.
\]
Since the maps $Y_n\hookrightarrow Y_{n+1}$ are cofibrations, so are the maps
$\mathrm{sp}_mY_n\hookrightarrow\mathrm{sp}_mY_{n+1}$;
for $n\geq m$ they are also weak equivalences
(because \(\mathrm{sp}_mY_{\bullet\geq m}\) represents \(\mathrm{core}_m(I_\bullet X_\bullet)_{\bullet\geq m}\),
whose transition maps are equivalences),
hence trivial cofibrations.
This concludes the proof because trivial cofibrations
are closed under sequential composites.
\end{proof}

We thus immediately obtain the main result of this article.

\begin{cor}
\label{LimitInd}
There is an equivalence
$\lim_{n\geq0}^{R} \mathscr Comp_n\simeq\mathscr Comp\noloc\vec{\mathrm{core}}$.
\end{cor}

\begin{proof}
  By definition of inductive tower.
\end{proof}

\begin{rmk}
\label{ComparisonWithFelixFramework}
  \Cref{LimitInd} can also be extracted from Loubaton's thesis~\cite{LoubatonThesis},
  by a detour through the $\Theta$-model
  for $(\infty,n)$-categories and $(\infty,\infty)$-categories.
  For the convenience of the reader, we give a short roadmap:
  \begin{itemize}
  \item
    \cite[Prop.~4.2.1.45]{LoubatonThesis} establishes that infinite $\Theta$-spaces
    form a model for the right limit of their finite stages.
  \item
    \cite[Prop.~3.4.3.9 and Thm.~3.4.3.14]{LoubatonThesis} show
    the Quillen adjunctions of \cite[Thm.~4.16]{ORquillen}
    (extended to infinity in \cite[Rem.~3.4.3.8]{LoubatonThesis})
    between the $\Theta$- and complicial models
    are Quillen equivalences at each finite stage and at infinity.
    The key ingredient in this comparison is a non-trivial combinatorial computation
    that guarantees that equivalences between complicial sets
    can be detected by globes \cite[Thm.~2.4.2.9]{LoubatonThesis}.
  \item
    The compatibility of these Quillen equivalences
    with the respective transition functors in the right limits
    (by passing to right adjoints in the compatibility with inclusions
    used in the proof of \cite[Lem.~3.4.3.13]{LoubatonThesis})
    allows a transfer of the claimed limiting statement.
  \end{itemize}
\end{rmk}

\subsection{Left limit: Truncation tower of a saturated complicial set}

\begin{prop}
\label{TruncationAdjunction}
For all $n\geq0$, there is an adjunction
\[
\adjunction{\tr_n}{\mathscr Comp}{\mathscr Comp_n}{I_n},
\]
where $\tr_n\coloneq \mathbb{L}\mathrm{th}_n$ denotes the $n$-th \emph{intelligent truncation},
which essentially enforces invertibility of all morphisms in dimensions greater than $n$.
\end{prop}

\begin{proof}
We have the Quillen pair
\[
\adjunction{\mathrm{th}_n}{m\sset_{\mathrm{cmp,sat}}}
  {m\sset_{(\infty,n)}}{\mathrm{sp}_n}
\]
from \cite[Prop.~4.5]{ORW2}.
Moreover, we have $\mathrm{sp}_n(X)=X$ for every $n$-complicial set,
hence $\mathbb{R}\mathrm{sp}_n\simeq I_n$, as desired.
\end{proof}

\begin{thm}
\label{VerificationDComplicial}
The family of adjunctions
\[
\{\adjunction{\tr_n}{\mathscr Comp}{\mathscr Comp_n}{I_n}\}_{n\geq0}
\]
is a dual preinductive tower for $\mathscr Comp$, but does not satisfy the dual of \ref{A4}.
\end{thm}

\begin{proof}
  Each $I_n\colon\mathscr Comp_n\to\mathscr Comp$ is fully faithful
  and admits the left adjoint $\tr_n$
  by \cref{TruncationAdjunction}.
  Moreover, the essential images are nested:
\[
I_n(\mathscr Comp_n)\subseteq I_{n+1}(\mathscr Comp_{n+1}),
\]
since every $n$-complicial set is also $(n+1)$-complicial. Thus $\{\mathscr Comp_n\}_{n\geq0}$ is a nested system of full reflective subcategories of $\mathscr Comp$, and the claim follows from the dual of \cref{prop:preinductive-reflective-subcats}.

It remains to show that the dual of \ref{A4} fails. Let $\cE$ be the walking coinductive equivalence from \cite[Constr.~1.5.13]{ORsurvey}, and consider the unique map
$f\colon N^{\mathrm{RS}}\cE\to\Delta[0]$.
By \cite[Thm.~4.21]{ORW2}, the marked simplicial set $N^{\mathrm{RS}}\cE$ is a saturated complicial set, the map $f$ is not a weak equivalence in $m\sset_{\mathrm{cmp,sat}}$, and, for every $n\geq0$, the map
$\mathrm{th}_nf\colon\mathrm{th}_nN^{\mathrm{RS}}\cE\to\Delta[0]$
is a weak equivalence in $m\sset_{\mathrm{cmp,sat}}$. Since $m\sset_{(\infty,n)}$ is a localization of $m\sset_{\mathrm{cmp,sat}}$, it follows that $\tr_n[f]$ is an equivalence in $\mathscr Comp_n$ for every $n\geq0$. However, $[f]$ is not an equivalence in $\mathscr Comp$. Hence the functors $\tr_n$ are not jointly conservative, so the dual of \ref{A4} does not hold.
\end{proof}

\begin{cor}
\label{LimitCoind}
There is an adjunction
$\adjunction{\vec{\mathrm{tr}}}{\mathscr Comp}{\lim_{n\geq0}^{L}\mathscr Comp_n}{}$,
which is a proper localization.
\end{cor}

\begin{proof}
  From \cref{VerificationDComplicial,LimitInd} we get the adjunction
  which is not an equivalence.
  Furthermore, we claim that the composite adjunction
  \begin{equation}
    \label{eq:composite-adjunction}
    {\lim^R_{n}\mathscr{C}omp_n}
    \rightleftarrows
    {\mathscr Comp}
    \rightleftarrows
    {\lim^L_n\mathscr Comp_n}
  \end{equation}
  is that of \cite[Proposition~4.6, Theorem~4.25]{ORW1}.
  One way to see this is to work in $\Pr^L$,
  where the right limit $\lim^R_n\mathscr Comp_n$ is simply the colimit along the inclusions,
  so that it suffices to compare the respective finite stage maps
  $\mathscr Comp_n \to \mathscr Comp_m$
  functorially in $n$ and $m$.
  By direct inspection one verifies that these finite stage maps
  are given by $L_mI_n$ in both adjunctions.

  Thus the claimed localization follows from \cite[Theorem~4.25]{ORW1}
  and \Cref{LimitInd}.
  One could have also used \cite[Theorem~1.10.4]{GH-hyper}
  instead, with a similar argument.
\end{proof}

\bibliographystyle{amsalpha}
\bibliography{ref}

@article {AraMaltsiniotisJoin,
    AUTHOR = {Ara, Dimitri and Maltsiniotis, Georges},
     TITLE = {Joint et tranches pour les $\infty$-cat\'{e}gories strictes},
   JOURNAL = {M\'{e}m. Soc. Math. Fr. (N.S.)},
  FJOURNAL = {M\'{e}moires de la Soci\'{e}t\'{e} Math\'{e}matique de France. Nouvelle S\'{e}rie},
    NUMBER = {165},
      YEAR = {2020},
     PAGES = {vi+213},
      ISSN = {0249-633-X},
      ISBN = {978-2-85629-921-0},
       DOI = {10.24033/msmf.473},
}

@article{CampionPasting,
      title={An $(\infty,n)$-categorical pasting theorem}, 
      author={Timothy Campion},
      year={2023},
      journal={\href{https://arxiv.org/abs/2311.00200v1}{arXiv:2311.00200v1}},
}

@article {henry-loubaton,
    AUTHOR = {Henry, Simon and Loubaton, F\'elix},
     TITLE = {An inductive model structure for strict {$\infty
              $}-categories},
   JOURNAL = {J. Pure Appl. Algebra},
  FJOURNAL = {Journal of Pure and Applied Algebra},
    VOLUME = {229},
      YEAR = {2025},
    NUMBER = {1},
     PAGES = {Paper No. 107859, 58},
      ISSN = {0022-4049,1873-1376},
   MRCLASS = {18N30 (18N40 18N65)},
  MRNUMBER = {4843919},
MRREVIEWER = {Julia\ Bergner},
       DOI = {10.1016/j.jpaa.2024.107859},
       URL = {https://doi.org/10.1016/j.jpaa.2024.107859},
}

@article{GepnerHeine,
      title={Oriented Category Theory}, 
      author={David Gepner and Hadrian Heine},
      year={2025},
      journal={\href{https://arXiv:2510.10504v1}{arXiv:2510.10504v1}},
}

@article{GH-hyper,
      title={Homotopy Posets, {P}ostnikov Towers, and Hypercompletions of $\infty$-Categories}, 
      author={David Gepner and Hadrian Heine},
      year={2026},
      journal={\href{https://arxiv.org/abs/2603.09903v2}{arXiv:2603.09903v2}},
}

@article {Goldthorpe,
    AUTHOR = {Goldthorpe, Zach},
     TITLE = {Homotopy theories of {$(\infty,\infty)$}-categories as
              universal fixed points with respect to weak enrichment},
   JOURNAL = {Int. Math. Res. Not. IMRN},
  FJOURNAL = {International Mathematics Research Notices. IMRN},
      YEAR = {2023},
    NUMBER = {22},
     PAGES = {19592--19640},
      ISSN = {1073-7928,1687-0247},
   MRCLASS = {18N65},
  MRNUMBER = {4669810},
MRREVIEWER = {Philippe\ Gaucher},
       DOI = {10.1093/imrn/rnad196},
       URL = {https://doi.org/10.1093/imrn/rnad196},
}

@article{JFR,
  title={{How to build a Hopf algebra}},
  author={Johnson-Freyd, Theo and Reutter, David},
  journal={\href{https://arXiv:2508.16787v1}{arXiv:2508.16787v1}},
  year={2025},
}

@book {LeinsterBook,
    AUTHOR = {Leinster, Tom},
     TITLE = {Higher operads, higher categories},
    SERIES = {London Mathematical Society Lecture Note Series},
    VOLUME = {298},
 PUBLISHER = {Cambridge University Press, Cambridge},
      YEAR = {2004},
     PAGES = {xiv+433},
      ISBN = {0-521-53215-9},
   MRCLASS = {18D50 (55P48)},
  MRNUMBER = {2094071},
MRREVIEWER = {R.\ H.\ Street},
       DOI = {10.1017/CBO9780511525896},
       URL = {https://doi.org/10.1017/CBO9780511525896},
}

@article{LMW,
  title={A folk model structure on omega-cat},
  author={Lafont, Yves and M{\'e}tayer, Fran{\c{c}}ois and Worytkiewicz, Krzysztof},
  journal={Advances in Mathematics},
  volume={224},
  number={3},
  pages={1183--1231},
  year={2010},
  publisher={Elsevier}
}

@article{LoubatonThesis,
  title={Theory and models of $(\infty,\omega)$-categories},
  author={Loubaton, F\'elix},
  journal={\href{https://arxiv.org/abs/2307.11931v2}{arXiv:2307.11931v2}},
  year={2023},
}

@article{Loubaton4,
  title={Categorical Theory of $(\infty,\omega)$-Categories},
  author={Loubaton, F\'elix},
  journal={\href{https://arxiv.org/abs/2406.05425}{arXiv:2406.05425v2}},
  year={2024},
}

@book {htt,
    AUTHOR = {Lurie, Jacob},
     TITLE = {Higher topos theory},
    SERIES = {Annals of Mathematics Studies},
    VOLUME = {170},
 PUBLISHER = {Princeton University Press, Princeton, NJ},
      YEAR = {2009},
     PAGES = {xviii+925},
      ISBN = {978-0-691-14049-0; 0-691-14049-9},
   MRCLASS = {18-02 (18B25 18E35 18G30 18G55 55U40)},
  MRNUMBER = {2522659},
MRREVIEWER = {Mark Hovey},
       DOI = {10.1515/9781400830558},
       URL = {http://dx.doi.org/10.1515/9781400830558},
}

@article {or,
    AUTHOR = {Ozornova, Viktoriya and Rovelli, Martina},
     TITLE = {Model structures for ({$\infty$},n)--categories on
              (pre)stratified simplicial sets and prestratified simplicial
              spaces},
   JOURNAL = {Algebr. Geom. Topol.},
  FJOURNAL = {Algebraic \& Geometric Topology},
    VOLUME = {20},
      YEAR = {2020},
    NUMBER = {3},
     PAGES = {1543--1600},
      ISSN = {1472-2747},
   MRCLASS = {18N65 (18N40 55U10 55U35)},
  MRNUMBER = {4105558},
       DOI = {10.2140/agt.2020.20.1543},
       URL = {https://doi.org/10.2140/agt.2020.20.1543},
}

@article {ORquillen,
    AUTHOR = {Ozornova, Viktoriya and Rovelli, Martina},
     TITLE = {A {Q}uillen adjunction between globular and complicial
              approaches to {$(\infty,n)$}-categories},
   JOURNAL = {Adv. Math.},
  FJOURNAL = {Advances in Mathematics},
    VOLUME = {421},
      YEAR = {2023},
     PAGES = {Paper No. 108980, 57},
      ISSN = {0001-8708,1090-2082},
   MRCLASS = {18N65 (18N50 55U10 55U35)},
  MRNUMBER = {4571188},
MRREVIEWER = {Philippe\ Gaucher},
       DOI = {10.1016/j.aim.2023.108980},
       URL = {https://doi.org/10.1016/j.aim.2023.108980},
}

@article{ORsurvey,
author = {Ozornova, Viktoriya and Rovelli, Martina},
title = {What is an equivalence in a higher category?},
journal = {Bulletin of the London Mathematical Society},
volume = {56},
number = {1},
pages = {1-58},
doi = {https://doi.org/10.1112/blms.12947},
url = {https://londmathsoc.onlinelibrary.wiley.com/doi/abs/10.1112/blms.12947},
eprint = {https://londmathsoc.onlinelibrary.wiley.com/doi/pdf/10.1112/blms.12947},
year = {2024}
}

@preamble{"\providecommand{\noopsort}[1]{}"}

@article{ORW1,
title={{\noopsort{ORW1}}{C}ores and localizations of $(\infty,\infty)$-categories},
author={Ozornova, Viktoriya and Rovelli, Martina and Walde, Tashi},
journal={\href{https://arxiv.org/abs/2603.11005v1}{arXiv:2603.11005v1}},
year={2026}
}

@article{ORW2,
title={{\noopsort{ORW2}}{A}bout the contractibility of the walking coinductive equivalence},
author={Ozornova, Viktoriya and Rovelli, Martina and Walde, Tashi},
journal={\href{https://arxiv.org/abs/2608.28355v1}{arXiv:2608.28355v1}},
year={2026}
}

@book {PolygraphBook,
    AUTHOR = {Ara, Dimitri and Burroni, Albert and Guiraud, Yves and Malbos,
              Philippe and M\'etayer, Fran{\c c}ois and Mimram, Samuel},
     TITLE = {Polygraphs: from rewriting to higher categories},
    SERIES = {London Mathematical Society Lecture Note Series},
    VOLUME = {495},
 PUBLISHER = {Cambridge University Press, Cambridge},
      YEAR = {2025},
     PAGES = {xx+648},
      ISBN = {978-1-009-49898-2},
   MRCLASS = {18-02 (18C35 18Gxx 18N30)},
  MRNUMBER = {4866320},
}

@incollection {EmilyNotes,
    AUTHOR = {Riehl, Emily},
     TITLE = {Complicial sets, an overture},
 BOOKTITLE = {2016 {MATRIX} annals},
    SERIES = {MATRIX Book Ser.},
    VOLUME = {1},
     PAGES = {49--76},
 PUBLISHER = {Springer, Cham},
      YEAR = {2018},
      ISBN = {978-3-319-72298-6; 978-3-319-72299-3},
   MRCLASS = {18G35},
  MRNUMBER = {3792516},
MRREVIEWER = {R.\ H.\ Street},
}

@book {RiehlVerityBook,
    AUTHOR = {Riehl, Emily and Verity, Dominic},
     TITLE = {Elements of {$\infty$}-category theory},
    SERIES = {Cambridge Studies in Advanced Mathematics},
    VOLUME = {194},
 PUBLISHER = {Cambridge University Press, Cambridge},
      YEAR = {2022},
     PAGES = {xix+759},
      ISBN = {978-1-108-83798-9},
   MRCLASS = {18-02 (18Nxx)},
  MRNUMBER = {4354541},
MRREVIEWER = {Philippe\ Gaucher},
       DOI = {10.1017/9781108936880},
       URL = {https://doi.org/10.1017/9781108936880},
}

@article{Stefanich,
  title={Presentable $(\infty, n)$-categories},
  author={Stefanich, Germ{\'a}n},
  journal={\href{https://arxiv.org/abs/2011.03035v1}{arXiv:2011.03035v1}},
  year={2020}
}

@article {StreetOrientedSimplexes,
    AUTHOR = {Street, Ross},
     TITLE = {The algebra of oriented simplexes},
   JOURNAL = {J. Pure Appl. Algebra},
  FJOURNAL = {Journal of Pure and Applied Algebra},
    VOLUME = {49},
      YEAR = {1987},
    NUMBER = {3},
     PAGES = {283--335},
      ISSN = {0022-4049},
   MRCLASS = {18G50 (18D05 18G55)},
  MRNUMBER = {920944},
MRREVIEWER = {Timothy Porter},
       DOI = {10.1016/0022-4049(87)90137-X},
       URL = {https://doi.org/10.1016/0022-4049(87)90137-X},
}

@misc{VeritySlides,
    AUTHOR = {Verity, Dominic},
    TITLE={A complicial compendium}, 
    year={2017},
    note={slides, available at \url{https://www.cirm-math.fr/ProgWeebly/Renc1773/Verity.pdf}, retrieved in October 2026},
}

\end{document}